\documentclass[12pt]{article}

\usepackage{amssymb,t1enc}
\newcommand{\wt}{\widetilde}
\newcommand{\R}{\mathbb R}
\newcommand{\C}{\mathbb C}
\newcommand{\Q}{\mathbb Q}
\newcommand{\A}{\textbf{\textit{A}}}
\newcommand{\B}{\textbf{\textit{B}}}
\newcommand{\F}{\textbf{\textit{F}}}
\newcommand{\K}{\textbf{\textit{K}}}
\newcommand{\sK}{\textbf{\textit{\scriptsize K}}}
\newcommand{\ii}{\textbf{\textit{i}}}
\begin{document}
\hfill
\vskip 4.0truecm
\thispagestyle{empty}
\par
\noindent
\centerline{{\large The Beckman-Quarles theorem for mappings from ${\C}^2$ to ${\C}^2$}}
\vskip 1.0truecm
\par
\noindent
\centerline{{\large Apoloniusz Tyszka}}
\vskip 1.0truecm
\par
\noindent
{\bf Abstract.} Let $\varphi:{\C}^2 \times {\C}^2 \to \C$,
$\varphi((x_1,x_2),(y_1,y_2))=(x_1-y_1)^2+(x_2-y_2)^2$.
We say that $f:{\C}^2 \to {\C}^2$ preserves distance $d \geq 0$,
if for each $X,Y \in {\C}^2$ $\varphi(X,Y)=d^2$ implies
$\varphi(f(X),f(Y))=d^2$. We prove that each unit-distance
preserving mapping $f:{\C}^2 \to {\C}^2$ has
a form $I \circ (\gamma,\gamma)$, where $\gamma: \C \to \C$ is
a field homomorphism and $I: {\C}^2 \to {\C}^2$ is
an affine mapping with orthogonal linear part.
We prove an analogous result for mappings from
${\K}^2$ to ${\K}^2$, where $\K$ is a commutative field such that
${\rm char}(\K) \not\in \{2,3,5\}$ and $-1$ is a square.
\vskip 1.0truecm
\par
The classical Beckman-Quarles theorem states that each
unit-distance preserving mapping from ${\R}^n$ to
${\R}^n$ ($n \geq 2$) is an isometry, see
\cite{Beckman-Quarles}--\cite{Everling}.
Let $\varphi:{\C}^2 \times {\C}^2 \to \C$,
$\varphi((x_1,x_2),(y_1,y_2))=(x_1-y_1)^2+(x_2-y_2)^2$.
We say that $f:{\C}^2 \to {\C}^2$ preserves distance~$d \geq 0$, if
for each $X,Y \in {\C}^2$ $\varphi(X,Y)=d^2$ implies $\varphi(f(X),f(Y))=d^2$.
If $f:{\C}^2 \to {\C}^2$ and for each $X,Y \in {\C}^2$
$\varphi(X,Y)=\varphi(f(X),f(Y))$, then $f$ is an affine mapping
with orthogonal linear part;
it follows from a general theorem proved in [3,~58~ff], see also [4, p. 30].
The author proved in \cite{Tyszka2}:
\vskip 0.2truecm
\par
\noindent
{\bf (1)}~~each unit-distance preserving mapping $f:{\C}^2 \to {\C}^2$
satisfies $\varphi(X,Y)=\varphi(f(X),f(Y))$
for all $X,Y \in {\C}^2$ with rational $\varphi(X,Y)$.
\newpage
\par
\noindent
{\bf Theorem 1.} If $f:{\C}^2 \to {\C}^2$ preserves unit distance,
$f((0,0))=(0,0)$, $f((1,0))=(1,0)$ and $f((0,1))=(0,1)$,
then there exists a field homomorphism $\rho: \R \to \C$
satisfying
$${\bf (2)}~~~~\forall x_1,x_2 \in \C~~f((x_1,x_2)) \in 
\{(\rho({\rm Re}(x_1))+\rho({\rm Im}(x_1)) \cdot \ii,~\rho({\rm Re}(x_2))+\rho({\rm Im}(x_2)) \cdot \ii),
$$
$$
~~~~~~~~~~~~~~~~~~~~~~~~~~~~~~~~~~~~~~~~~~~~~~~(\rho({\rm Re}(x_1))-\rho({\rm Im}(x_1)) \cdot \ii,~\rho({\rm Re}(x_2))-\rho({\rm Im}(x_2)) \cdot \ii)\}.
$$
\par
\noindent
{\it Proof.} Obviously, $g=f_{|{\R}^2}:{\R}^2 \to {\C}^2$
preserves unit distance. The author proved in~\cite{Tyszka1}
that such a $g$ has a form $I \circ (\rho,\rho)$, where
$\rho: \R \to \C$ is a field homomorphism and
$I: {\C}^2 \to {\C}^2$ is an affine mapping
with orthogonal linear part. Since $f((0,0))=(0,0)$,
$f((1,0))=(1,0)$, $f((0,1))=(0,1)$,
we conclude that $f_{|{\R}^2}=(\rho,\rho)$.
From this, condition~{\bf (2)} holds true
if $(x_1,x_2) \in {\R}^2$. Assume now that
$(x_1,x_2) \in {\C}^2 \setminus {\R}^2$.
Let $x_1=a_1+b_1 \cdot \ii$, $x_2=a_2+b_2 \cdot \ii$,
where $a_1,b_1,a_2,b_2 \in \R$, and, for example $b_1 \neq 0$.
\par
\noindent
For each $t \in \R$
\par
\centerline{$\varphi((a_1+b_1 \cdot \ii,~a_2+b_2 \cdot \ii),~(a_1+tb_2,~a_2-tb_1))=(t^2-1)(b_1^2+b_2^2)$.}
\par
\noindent
By this and {\bf (1)}:
\vskip 0.2truecm
\par
\noindent
{\bf (3)}~~for each $t \in \R$ with rational $(t^2-1)(b_1^2+b_2^2)$ we have:
\par
\centerline{$\varphi(f((a_1+b_1 \cdot \ii,~a_2+b_2 \cdot \ii)),~f((a_1+tb_2,~a_2-tb_1)))=(t^2-1)(b_1^2+b_2^2)$.}
\vskip 0.2truecm
\par
\noindent
Let $f((a_1+b_1 \cdot \ii,~a_2+b_2 \cdot \ii))=(y_1,y_2)$.
From {\bf (3)} and $f_{|{\R}^2}=(\rho,\rho)$ we obtain:
\vskip 0.2truecm
\par
\noindent
{\bf (4)}~~for each $t \in \R$ with rational $(t^2-1)(b_1^2+b_2^2)$ we have:
\par
~~~~$(y_1-\rho(a_1)-\rho(t)\rho(b_2))^2+(y_2-\rho(a_2)+\rho(t)\rho(b_1))^2=(t^2-1)(b_1^2+b_2^2)$.
\vskip 0.2truecm
\par
\noindent
For each $t \in \R$ with rational $(t^2-1)(b_1^2+b_2^2)$ we have:
\vskip 0.2truecm
\par
\noindent
\centerline{$(t^2-1)(b_1^2+b_2^2)=\rho((t^2-1)(b_1^2+b_2^2))=
(\rho(t)^2-1)(\rho(b_1)^2+\rho(b_2)^2)$.}
\vskip 0.2truecm
\par
\noindent
By this and {\bf (4)}:
\vskip 0.2truecm
\par
\noindent
{\bf (5)}~~for each $t\in \R$ with rational $(t^2-1)(b_1^2+b_2^2)$ we have:
\par
\footnotesize
~~$(y_1-\rho(a_1))^2+(y_2-\rho(a_2))^2+\rho(b_1)^2+\rho(b_2)^2~+
2\rho(t) \cdot (\rho(b_1)(y_2-\rho(a_2))-\rho(b_2)(y_1-\rho(a_1)))=0$.
\normalsize
\vskip 0.2truecm
\par
\noindent
There are infinitely many $t \in \R$
with rational $(t^2-1)(b_1^2+b_2^2)$ and $\rho$ is injective.
From these two facts and {\bf (5)}, we obatin:
\vskip 0.2truecm
\par
\noindent
{\bf (6)}~~$\rho(b_1)(y_2-\rho(a_2))-\rho(b_2)(y_1-\rho(a_1))=0$
\vskip 0.2truecm
\par
\noindent
and
\vskip 0.2truecm
\par
\noindent
{\bf (7)}~~$(y_1-\rho(a_1))^2+(y_2-\rho(a_2))^2+\rho(b_1)^2+\rho(b_2)^2=0$.
\newpage
\par
\noindent
By {\bf (6)}:
\vskip 0.2truecm
\par
\noindent
{\bf (8)}~~$y_2-\rho(a_2)=\frac{\textstyle \rho(b_2)}{\textstyle \rho(b_1)} \cdot (y_1-\rho(a_1))$.
\vskip 0.2truecm
\par
\noindent
Applying {\bf (8)} to {\bf (7)} we get:
$$
(y_1-\rho(a_1))^2+\frac{\rho(b_2)^2}{\rho(b_1)^2} \cdot (y_1-\rho(a_1))^2 + \rho(b_1)^2+\rho(b_2)^2=0.
$$
It gives $\left(\frac{\textstyle (y_1-\rho(a_1))^2}{\textstyle \rho(b_1)^2}+1 \right) \cdot (\rho(b_1)^2+\rho(b_2)^2)=0$.
Since $\rho(b_1)^2+\rho(b_2)^2 \neq 0$, we get
$$\underbrace{y_1=\rho(a_1)+\rho(b_1) \cdot \ii}_{\rm case~1}~~{\rm or}~~\underbrace{y_1=\rho(a_1)-\rho(b_1) \cdot \ii}_{\rm case~2}.$$
\par
\noindent
In case 1, by {\bf (8)}
\vskip 0.2truecm
\par
\noindent
\small
$y_2=\rho(a_2)+\frac{\textstyle \rho(b_2)}{\textstyle \rho(b_1)} \cdot (y_1-\rho(a_1))=\rho(a_2)+\frac{\textstyle \rho(b_2)}{\textstyle \rho(b_1)} \cdot (\rho(a_1)+\rho(b_1) \cdot \ii -\rho(a_1))=\rho(a_2)+\rho(b_2) \cdot \ii$.
\normalsize
\vskip 0.2truecm
\par
\noindent
In case 2, by {\bf (8)}
\vskip 0.2truecm
\small
\par
\noindent
$y_2=\rho(a_2)+\frac{\textstyle \rho(b_2)}{\textstyle \rho(b_1)} \cdot (y_1-\rho(a_1))=\rho(a_2)+\frac{\textstyle \rho(b_2)}{\textstyle \rho(b_1)} \cdot (\rho(a_1)-\rho(b_1) \cdot \ii -\rho(a_1))=\rho(a_2)-\rho(b_2) \cdot \ii$.
\normalsize
\vskip 0.2truecm
\par
\noindent
The proof is completed.
\vskip 0.2truecm
Let
$f:{\C}^2 \to {\C}^2$ preserves unit distance,
$f((0,0))=(0,0)$, $f((1,0))=(1,0)$ and $f((0,1))=(0,1)$.
Theorem 1 provides a field homomorphism $\rho:\R \to \C$
satisfying~{\bf (2)}.
By Theorem 1 the sets
$$
\A=\{(x_1,x_2) \in {\C}^2:
f((x_1,x_2))=(\rho({\rm Re}(x_1))+\rho({\rm Im}(x_1)) \cdot \ii,~\rho({\rm Re}(x_2))+\rho({\rm Im}(x_2)) \cdot \ii)\}
$$
and
$$
\B=\{(x_1,x_2) \in {\C}^2:
f((x_1,x_2))=(\rho({\rm Re}(x_1))-\rho({\rm Im}(x_1)) \cdot \ii,~\rho({\rm Re}(x_2))-\rho({\rm Im}(x_2)) \cdot \ii)\}
$$
satisfy $\A \cup \B={\C}^2$. The mapping
$$
\C \ni x \stackrel{\theta}{\longrightarrow} \rho({\rm Re}(x))+\rho({\rm Im}(x)) \cdot \ii \in \C
$$
is a field homomorphism, $\theta$ extends $\rho$,
$$
\A=\{(x_1,x_2) \in {\C}^2: f((x_1,x_2))=(\theta(x_1),\theta(x_2))\}.
$$
\par
\noindent
The mapping
$$
\C \ni x \stackrel{\zeta}{\longrightarrow} \rho({\rm Re}(x))-\rho({\rm Im}(x)) \cdot \ii \in \C
$$
is a field homomorphism, $\zeta$ extends $\rho$,
$$
\B=\{(x_1,x_2) \in {\C}^2: f((x_1,x_2))=(\zeta(x_1),\zeta(x_2))\}.
$$
We would like to prove $f=(\theta,\theta)$ or $f=(\zeta,\zeta)$;
we will prove it later in Theorem 2.
\vskip 0.2truecm
\par
\noindent
Let $\psi:{\C}^2 \times {\C}^2 \to \R$,
$\psi((x_1,x_2),(y_1,y_2))={\rm Im}(x_1) \cdot {\rm Im}(y_1)+{\rm Im}(x_2) \cdot {\rm Im}(y_2)$.
\vskip 0.2truecm
\par
\noindent
{\bf Lemma 1.} If $x_1,x_2,y_1,y_2 \in \C$, $\varphi((x_1,x_2),(y_1,y_2)) \in \Q$
and $\psi((x_1,x_2),(y_1,y_2)) \neq 0$, then
\par
\noindent
${\bf (9)}~~~~~~~~~~~~~~~~~~~~~~~~~~~~~(y_1,y_2) \in \A~~{\rm implies}~~(x_1,x_2) \in \A$
\par
\noindent
and
\par
\noindent
${\bf (10)}~~~~~~~~~~~~~~~~~~~~~~~~~~~~(y_1,y_2) \in \B~~{\rm implies}~~(x_1,x_2) \in \B$.
\vskip 0.2truecm
\par
\noindent
{\it Proof.} We prove only {\bf (9)}, the proof of {\bf (10)}
follows analogically.
Let $\varphi((x_1,x_2),(y_1,y_2))=r \in \Q$.
Assume, on the contrary, that $(y_1,y_2) \in \A$ and
$(x_1,x_2) \not\in \A$. Since $\A \cup \B={\C}^2$,
$(x_1,x_2) \in \B$.
Let $x_1=a_1+{b_1} \cdot \ii$,~~$x_2=a_2+{b_2} \cdot \ii$,
~~$y_1=\wt{a_1}+\wt{b_1} \cdot \ii$,~~$y_2=\wt{a_2}+\wt{b_2} \cdot \ii$,
where $a_1,b_1,a_2,b_2,\wt{a_1},\wt{b_1},\wt{a_2},\wt{b_2} \in \R$.
By {\bf (1)}:
\vskip 0.2truecm
\par
\noindent
\centerline{$r=\varphi((x_1,x_2),(y_1,y_2))=\varphi(f((x_1,x_2)),f((y_1,y_2)))=$}
\par
\noindent
{\bf (11)}
\par
\noindent
\centerline{$(\rho(a_1)-\rho(b_1) \cdot \ii-\rho(\wt{a_1})-\rho(\wt{b_1}) \cdot \ii)^2
+
(\rho(a_2)-\rho(b_2) \cdot \ii-\rho(\wt{a_2})-\rho(\wt{b_2}) \cdot \ii)^2$.}
\vskip 0.3truecm
\par
\noindent
Since $r \in \Q$,
\par
\noindent
\centerline{
$r=\theta(r)=\theta((a_1+b_1 \cdot \ii - \wt{a_1} - \wt{b_1} \cdot \ii)^2+
(a_2+b_2 \cdot \ii - \wt{a_2} - \wt{b_2} \cdot \ii)^2)=$}
\par
\noindent
{\bf (12)}
\par
\noindent
$~~~~~~~~~~~~(\rho(a_1)+\rho(b_1) \cdot \ii-\rho(\wt{a_1})-\rho(\wt{b_1}) \cdot \ii)^2+
(\rho(a_2)+\rho(b_2) \cdot \ii-\rho(\wt{a_2})-\rho(\wt{b_2}) \cdot \ii)^2$.
\vskip 0.2truecm
\par
\noindent
Subtracting {\bf (11)} and {\bf (12)} by sides we obtain:
$$
2\rho(b_1) \cdot \ii \cdot 
(2\rho(\wt{b_1}) \cdot \ii-2\rho(a_1)+2\rho(\wt{a_1}))+
2\rho(b_2) \cdot \ii \cdot 
(2\rho(\wt{b_2}) \cdot \ii-2\rho(a_2)+2\rho(\wt{a_2}))=0.
$$
\vskip 0.2truecm
\par
\noindent
Thus
\vskip 0.2truecm
\par
\noindent
${\bf (13)}~~~~~~~~~~~~~~~~~~~
-\rho(b_1\wt{b_1}+b_2\wt{b_2})=
\rho(b_1(a_1-\wt{a_1})+b_2(a_2-\wt{a_2})) \cdot \ii$.
\vskip 0.2truecm
\par
\noindent
Squaring both sides of {\bf (13)} we get:
\par
\noindent
\vskip 0.3truecm
\centerline{
$\rho((b_1\wt{b_1}+b_2\wt{b_2})^2+(b_1(a_1-\wt{a_1})+b_2(a_2-\wt{a_2}))^2)=0$,}
\vskip 0.3truecm
\par
\noindent
so in particular $\psi((x_1,x_2),(y_1,y_2))=b_1\wt{b_1}+b_2\wt{b_2}=0$, a contradiction.
\vskip 0.2truecm
\par
\noindent
The next lemma is obvious.
\vskip 0.2truecm
\par
\noindent
{\bf Lemma 2.} For each $S,T \in {\R}^2$ there exist $n \in \{1,2,3,...\}$
and $P_1,...,P_n \in {\R}^2$ such that
$||S-P_1||=||P_1-P_2||=...=||P_{n-1}-P_n||=||P_n-T||=1$.
\vskip 0.2truecm
\par
\noindent
{\bf Lemma 3.} For each $X \in {\C}^2 \setminus {\R}^2$
\vskip 0.2truecm
\par
\noindent
\centerline{$(\ii,\ii) \in \A$ implies $X \in \A$}
\par
\noindent
and
\par
\noindent
\centerline{$(\ii,\ii) \in \B$ implies $X \in \B$.}
\vskip 0.2truecm
\par
\noindent
{\it Proof.} Let $X=(a_1+b_1 \cdot \ii,~a_2+b_2 \cdot \ii)$, where
$a_1,b_1,a_2,b_2 \in \R$. Since $X \in {\C}^2 \setminus {\R}^2$,
$b_1 \neq 0$ or $b_2 \neq 0$. Assume that $b_1 \neq 0$, when
$b_2 \neq 0$ the proof is analogous.
The points $S=\left(a_1+\sqrt{1+b_2^2},~a_2+\sqrt{1+(b_1-1)^2}\right)$ and
$T=\left(\sqrt{2},0\right)$ belong to ${\R}^2$. Applying Lemma 2 we find
$P_1,...,P_n \in {\R}^2$ satisfying
$||S-P_1||=||P_1-P_2||=...=||P_{n-1}-P_n||=||P_n-T||=1$.
The points
\vskip 0.08truecm
\par
\noindent
$X_1=X$,
\vskip 0.08truecm
\par
\noindent
$X_2=\left(a_1+\sqrt{1+b_2^2}+b_1 \cdot \ii,~a_2\right)$,
\vskip 0.08truecm
\par
\noindent
$X_3=S+(\ii,0)=\left(a_1+\sqrt{1+b_2^2}+\ii,~a_2+\sqrt{1+(b_1-1)^2}\right)$,
\vskip 0.08truecm
\par
\noindent
$X_4=P_1+(\ii,0)$,
\vskip 0.08truecm
\par
\noindent
$X_5=P_2+(\ii,0)$,
\vskip 0.08truecm
\par
\noindent
. . . . . . . . . . .
\vskip 0.08truecm
\par
\noindent
$X_{n+3}=P_n+(\ii,0)$,
\vskip 0.08truecm
\par
\noindent
$X_{n+4}=T+(\ii,0)=\left(\sqrt{2}+\ii,~0\right)$,
\vskip 0.08truecm
\par
\noindent
$X_{n+5}=(\ii,\ii)$
\vskip 0.08truecm
\par
\noindent
satisfy:
\vskip 0.08truecm
\par
\noindent
for each $k \in \{2,3,...,n+5\}$ $\varphi(X_{k-1},X_k)=1$,
\par
\noindent
$\psi(X_1,X_2)=b_1^2 \neq 0$, $\psi(X_2,X_3)=b_1 \neq 0$, for 
each $k \in \{4,5,...,n+5\}$ $\psi(X_{k-1},X_k)=1$.
\vskip 0.2truecm
\par
\noindent
By Lemma 1 for each $k \in \{2,3,...,n+5\}$
\vskip 0.2truecm
\par
\centerline{$X_k \in \A$ implies $X_{k-1} \in \A$}
\par
\noindent
and
\par
\noindent
\centerline{$X_k \in \B$ implies $X_{k-1} \in \B$.}
\vskip 0.2truecm
\par
\noindent
Therefore, $(\ii,\ii)=X_{n+5} \in \A$ implies $X=X_1 \in \A$, and also, 
$(\ii,\ii)=X_{n+5} \in \B$ implies $X=X_1 \in \B$.
\newpage
\par
\noindent
{\bf Theorem 2.} If $f:{\C}^2 \to {\C}^2$ preserves unit distance,
$f((0,0))=(0,0)$, $f((1,0))=(1,0)$ and $f((0,1))=(0,1)$,
then there exists a field homomorphism $\gamma: \C \to \C$
satisfying $f=(\gamma,\gamma)$.
\vskip 0.2truecm
\par
\noindent
{\it Proof.} By Lemma 3
\par
\noindent
\centerline{$(\ii,\ii) \in \A$~~implies~~${\C}^2 \setminus {\R}^2 \subseteq \A$}
\par
\noindent
and
\par
\noindent
\centerline{$(\ii,\ii) \in \B$~~implies~~${\C}^2 \setminus {\R}^2 \subseteq \B$.}
\vskip 0.2truecm
\par
\noindent
Obviously, ${\R}^2 \subseteq \A$ and ${\R}^2 \subseteq \B$. Therefore,
\par
\noindent
\centerline{$\A={\C}^2$ and $f=(\theta,\theta)$,~~if $(\ii,\ii) \in \A$,}
\par
\noindent
and also,
\par
\noindent
\centerline{$\B={\C}^2$ and $f=(\zeta,\zeta)$,~~if $(\ii,\ii) \in \B$.}
\vskip 0.6truecm
\par
\noindent
As a corollary of Theorem 2 we get:
\vskip 0.2truecm
\par
\noindent
{\bf Theorem 3.} Each unit-distance preserving mapping
$f:{\C}^2 \to {\C}^2$ has a form $I \circ (\gamma,\gamma)$,
where $\gamma: \C \to \C$ is a field homomorphism and
$I:{\C}^2 \to {\C}^2$ is an affine mapping with orthogonal linear part.
\vskip 0.20truecm
\par
\noindent
{\it Proof.} By {\bf (1)}:
\par
\noindent
\centerline{$1=\varphi((0,0),(1,0))=\varphi(f((0,0)),f((1,0)))$,}
\par
\noindent
\centerline{$1=\varphi((0,0),(0,1))=\varphi(f((0,0)),f((0,1)))$,}
\par
\noindent
\centerline{$2=\varphi((1,0),(0,1))=\varphi(f((1,0)),f((0,1)))$.}
\vskip 0.20truecm
\par
\noindent
By the above equalities there exists an
affine mapping $J:{\C}^2 \to {\C}^2$
with orthogonal linear part such that
$J(f((0,0)))=(0,0)$, $J(f((1,0)))=(1,0)$,
$J(f((0,1)))=(0,1)$. By Theorem 2 there exists
a field homomorphism $\gamma: \C \to \C$ satisfying $J \circ f=(\gamma,\gamma)$,
so $f=J^{-1} \circ (\gamma,\gamma)$.
\vskip 0.20truecm
\par
\noindent
Obviously, Theorem 3 implies {\bf (1)}. The author proved in \cite{Tyszka3}:
\vskip 0.20truecm
\par
\noindent
{\bf (14)}~~if $n \geq 2$ and a continuous $f:{\C}^n \to {\C}^n$ 
preserves unit distance, then $f$ has a
form $I \circ (\rho,...,\rho)$, where $I:{\C}^n \to {\C}^n$ is an affine
mapping with orthogonal linear
part and $\rho:\C \to \C$ is the identity
or the complex conjugation.
\vskip 0.2truecm
\par
\noindent
The only continuous endomorphisms of $\C$ are the identity and the complex
conjugation, see [6, Lemma 1, p. 356]. Therefore, Theorem 3 implies 
{\bf (14)} restricted to $n=2$.
\newpage
Let $\K$ be a commutative field, ${\rm char}(\K) \not\in \{2,3,5\}$.
Let $d:{\K}^2 \times {\K}^2 \to \K$ denote the
Lorentz-Minkowski distance defined by
$d((x_1,x_2),(y_1,y_2))=(x_1-y_1) \cdot (x_2-y_2)$.
H. Schaeffer proved in [7, Satz 1, Satz 2, Satz 3]:
\vskip 0.20truecm
\par
\noindent
{\bf (15)}~~if $f:{\K}^2 \to {\K}^2$ preserves the Lorentz-Minkowski
distance~$1$, $f((0,0))=(0,0)$ and $f((1,1))=(1,1)$, then there exists
a field homomorphism $\sigma: \K \to \K$ satisfying
$\forall x_1,x_2 \in \K~f((x_1,x_2))=(\sigma(x_1),\sigma(x_2))~
{\rm or}~\forall x_1,x_2 \in \K~f((x_1,x_2))=(\sigma(x_2),\sigma(x_1))$.
\normalsize
\vskip 0.20truecm
\par
\noindent
Unfortunately, the proof of Satz 3 in \cite{Schaeffer} is complicated,
the main part of this proof was constructed using computer software.
\vskip 0.20truecm
\par
Let $\varphi_{\sK}:{\K}^2 \times {\K}^2 \to \K$,
$\varphi_{\sK}((x_1,x_2),(y_1,y_2))=(x_1-y_1)^2+(x_2-y_2)^2$.
Theorem~4 generalizes Theorem 3.
\vskip 0.20truecm
\par
\noindent
{\bf Theorem 4.} Let there exists $i \in \K$ such that $i^2+1=0$.
Let $f:{\K}^2 \to {\K}^2$ preserves unit distance defined by
$\varphi_{\sK}$. We claim that $f$ has a form $I \circ (\sigma, \sigma)$,
where $\sigma: \K \to \K$ is a field homomorphism and $I:{\K}^2 \to {\K}^2$
is an affine mapping with orthogonal linear part.
\vskip 0.2truecm
\par
\noindent
{\it Proof.} Assume that $f((0,0))=(0,0)$. The mappings
\vskip 0.3truecm
\par
\noindent
\centerline{${\K}^2 \ni (x_1,x_2) \stackrel{\xi}{\longrightarrow} (x_1+i \cdot x_2,~x_1-i \cdot x_2) \in {\K}^2$}
\par
\noindent
and
\par
\noindent
\centerline{~~~${\K}^2 \ni (x_1,x_2) \stackrel{\eta}{\longrightarrow} \left(\frac{1}{2}x_1+\frac{1}{2}x_2,~-\frac{i}{2}x_1+\frac{i}{2}x_2 \right) \in {\K}^2$}
\vskip 0.3truecm
\par
\noindent
satisfy:
\par
\noindent
~~~~~~~~~~$\eta \circ \xi=\xi \circ \eta={\rm id}({\K}^2)$,
\par
\noindent
~~~~~~~~~~$\forall x_1,x_2,y_1,y_2 \in \K~~\varphi_{\sK}((x_1,x_2),(y_1,y_2))=d(\xi((x_1,x_2)),\xi((y_1,y_2)))$,
\par
\noindent
~~~~~~~~~~$\forall x_1,x_2,y_1,y_2 \in \K~~d((x_1,x_2),(y_1,y_2))=\varphi_{\sK}(\eta((x_1,x_2)),\eta((y_1,y_2)))$.
\vskip 0.3truecm
\par
\noindent
Therefore, $\xi \circ f \circ \eta: {\K}^2 \to {\K}^2$
preserves the Lorentz-Minkowski distance~$1$.
Obviously, $(\xi \circ f \circ \eta)((0,0))=(0,0)$.
Let $(\xi \circ f \circ \eta)((1,1))=(a,b) \in {\K}^2$. We have:
$1=d((1,1),(0,0))=d((\xi \circ f \circ \eta)((1,1)),(\xi \circ f \circ \eta)((0,0)))=d((a,b),(0,0))=a \cdot b$.
Hence $b=\frac{1}{a}$.
For each $z \in \K \setminus \{0\}$ the mapping
\vskip 0.3truecm
\par
\centerline{${\K}^2 \ni (x,y) \stackrel{\lambda(z)}{\longrightarrow} (\frac{x}{z}, z \cdot y) \in {\K}^2$}
\vskip 0.3truecm
\par
\noindent
preserves all Lorentz-Minkowski distances,
$\lambda(\frac{1}{z}) \circ \lambda(z)=\lambda(z) \circ \lambda(\frac{1}{z})={\rm id}({\K}^2)$.
The mapping $\lambda(a) \circ \xi \circ f \circ \eta:{\K}^2 \to {\K}^2$
preserves the Lorentz-Minkowski distance~$1$,
$(\lambda(a) \circ \xi \circ f \circ \eta)((0,0))=(0,0)$ and
$(\lambda(a) \circ \xi \circ f \circ \eta)((1,1))=(1,1)$.
By {\bf (15)} there exists a field homomorphism $\sigma: \K \to \K$ satisfying
\vskip 0.3truecm
\par
\noindent
\centerline{$
\underbrace{\lambda(a) \circ \xi \circ f \circ \eta=(\sigma,\sigma)}_{\rm case~1}
~~~~{\rm or}~~~~
\underbrace{\lambda(a) \circ \xi \circ f \circ \eta=h \circ (\sigma,\sigma)}_{\rm case~2}~,
$}
\vskip 0.3truecm
\par
\noindent
where $h:{\K}^2 \to {\K}^2$, $h((x_1,x_2))=(x_2,x_1)$.
\vskip 0.3truecm
\par
\noindent
In case~1:~~$f=\eta \circ \lambda(\frac{1}{a}) \circ (\sigma,\sigma) \circ \xi=f_1 \circ (\sigma,\sigma)$,
where $f_1:{\K}^2 \to {\K}^2$,
\vskip 0.3truecm
\par
\noindent
\centerline{
$f_1((x_1,x_2))=
\left(\left(\frac{a}{2}+\frac{1}{2a}\right)\cdot x_1+\left(\frac{a}{2}-\frac{1}{2a}\right) \sigma(i) \cdot x_2,
~~
-\left(\frac{a}{2}-\frac{1}{2a}\right)i \cdot x_1-\left(\frac{a}{2}+\frac{1}{2a}\right)i\sigma(i) \cdot x_2 \right)
$.}
\vskip 0.3truecm
\par
\noindent
In case~2:~~$f=\eta \circ \lambda(\frac{1}{a}) \circ h \circ (\sigma,\sigma) \circ \xi=f_2 \circ (\sigma, \sigma)$,
where $f_2:{\K}^2 \to {\K}^2$,
\vskip 0.3truecm
\par
\noindent
\centerline{
$f_2((x_1,x_2))=
\left(\left(\frac{a}{2}+\frac{1}{2a}\right)\cdot x_1 - \left(\frac{a}{2}-\frac{1}{2a}\right) \sigma(i) \cdot x_2,
~~
-\left(\frac{a}{2}-\frac{1}{2a}\right)i \cdot x_1+\left(\frac{a}{2}+\frac{1}{2a}\right)i\sigma(i) \cdot x_2 \right)
$.}
\vskip 0.3truecm
\par
\noindent
The mappings $f_1$ and $f_2$ are linear and orthogonal. The proof is completed.

\par
\noindent
Apoloniusz Tyszka\\
Technical Faculty\\
Hugo Ko\l{}\l{}\k{a}taj University\\
Balicka 104, 30-149 Krak\'ow, Poland\\
E-mail address: {\it rttyszka@cyf-kr.edu.pl}
\end{document}